\documentclass[12pt]{amsart}%
\usepackage{cite}
\usepackage{graphicx}
\usepackage{amscd}
\usepackage{amsmath}
\usepackage{amsfonts}
\usepackage{amssymb}%
\usepackage{enumitem}
\usepackage{url}
\providecommand{\U}[1]{\protect\rule{.1in}{.1in}}
\newtheorem{theorem}{Theorem}
\theoremstyle{plain}

\newtheorem{definition}{Definition}[section]

\theoremstyle{definition}
\newtheorem{example}{Example}[section]

\newtheorem{remark}{Remark}[section]

\numberwithin{equation}{section}
\numberwithin{theorem}{section}

\DeclareMathOperator{\arccot}{arccot}
\begin{document}
\title{Carl St\o rmer and his Numbers} 
\author{Matthew Kroesche}
\address{Department of Mathematics, Sid Richardson Building, 1410 S.4th Street, Waco,
TX 76706.}
\email{matthew\_kroesche@baylor.edu}
\author{Lance L. Littlejohn}
\address{Department of Mathematics, Sid Richardson Building, 1410 S.4th Street, Waco,
TX 76706.}
\email{lance\_littlejohn@baylor.edu}
\author{Graeme Reinhart}
\address{Department of Mathematics, Sid Richardson Building, 1410 S.4th Street, Waco,
TX 76706.}
\email{graeme\_reinhart1@baylor.edu}

\begin{abstract}
In many proofs of Fermat's Two Squares Theorem, the smaller of the two least residue solutions $x_0$ of the quadratic congruence $x^2 \equiv -1 \bmod p$ plays an essential role; here $p$ is prime and $p \equiv 1 \bmod 4$. Such an $x_0$ is called a St\o rmer number, named after the Norwegian mathematician and astronomer Carl St\o rmer (1874-1957). In this paper, we establish necessary and sufficient conditions for $x_0 \in \mathbb{N}$ to be a St\o rmer number of some prime $p \equiv 1 \bmod 4$. St\o rmer's main interest in his investigations of St\o rmer numbers stemmed from his study of identities expressing $\pi$ as finite linear combinations of certain values of the Gregory-MacLaurin series for $\arctan(1/x)$. Since less than 600 digits of $\pi$ were known by 1900,  approximating $\pi$ was an important topic. One such identity, discovered by St\o rmer in 1896, was used by Yasumasa Kanada and his team in 2002 to obtain 1.24 trillion digits of $\pi$. We also provide a new proof of Størmer's theorem on the decomposition of Gregory numbers, and discuss this theorem's connection to approximations of $\pi$.  
\end{abstract}
\date{}
\maketitle

\section{Introduction and Motivation}

Fermat's Two Squares Theorem, which G. H. Hardy said is\textit{ ``ranked, very justly, as one of the finest of arithmetic''} (see \cite{HardyApology}) states
\begin{theorem}
A prime number $p$ can be written as the sum of squares of two positive
integers only when $p\equiv 1\bmod 4.$ Furthermore, in this case, there exists
unique positive integers $a$ and $b$ such that%
\[
p=a^{2}+b^{2}.
\]
\end{theorem}

From a historical point of view, L. E. Dickson \cite{CELV-MAA-1999} attributes the statement of the above theorem
to A. Girard in 1632. In a letter to Mersenne, Fermat's announced the proof in
1640; E. T. Bell \cite[p. 89]{Bell} and Dickson \cite[p. 228]{Dickson} say
that Fermat had an irrefutable proof by the method of infinite descent.
However, no record of Fermat's proof has been found. In 1749, Euler gave the
first recorded proof of the existence part of Fermat's theorem. Bell further
writes \textit{``It was first proved by the great Euler in 1749 after
he struggled, off and on, for seven years to find a proof''}.
The first proof of uniqueness of the positive integers was given by Gauss in 1801.

Most proofs of the existence of the two squares initially involve a study of the quadratic
congruence $x^{2}\equiv-1 \bmod p$. For example, Hermite \cite{Hermite} and Serret \cite{Serret} who, independently in 1848, gave proofs of Fermat's Two Squares Theorem assuming the least residue solution $x_0$ of $x^2 \equiv -1 \bmod p$ is known. Then, using the Euclidean Algorithm, they produced algorithms to find the positive integers $a,b$ satisfying Fermat's Two Squares Theorem. For further information, see the contribution by Brillhart \cite{Brillhart1972} who gave an account of Serret's and Hermite's algorithms and improved both. 

In 1855, British mathematician Henry John Stephen Smith gave an elegant existence proof
of Fermat's Two Squares Theorem using the Euclidean Algorithm and elementary
theory of determinants and continuants; see \cite{CELV-MAA-1999} and \cite{HJS Smith 1855}. We briefly describe Smith's method which also emphasizes the importance of the quadratic congruence $x^2 \equiv-1\bmod p$. For explicit details of Smith's proof, see \cite{CELV-MAA-1999} and \cite{HJS Smith 1855}. 

For $m \in \mathbb{N}$ and positive integers $q_1, q_2, \ldots, q_m \in \mathbb{N}$, the \textit{continuant} $[q_1,q_2,\ldots,q_m]$ of length $m$ is defined by $[q_{1}]=q_{1}$ and, for $m>1$, it is the determinant
\begin{equation}
\lbrack q_{1},q_{2},\ldots,q_{m}]=\left\vert
\begin{array}
[c]{rrrrrr}%
q_{1} & 1 & 0 & \cdots & 0 & 0\\
-1 & q_{2} & 1 & 0 & \cdots & 0\\
0 & -1 & q_{3} & 1 & \cdots & 0\\
\vdots &  &  &  &  & \vdots\\
0 & 0 & 0 & \cdots & q_{m-1} & 1\\
0 & 0 & 0 & \cdots & -1 & q_{m}%
\end{array}
\right\vert . \nonumber%
\end{equation}
It is not difficult to see that
\[
[q_1,q_2,\ldots, q_{m-1},q_m]=[q_m,q_{m-1},\ldots,q_2,q_1]; 
\]
for additional properties of continuants, see \cite{CELV-MAA-1999} and \cite{Muir and Metzler}. 
Continuants arise naturally from applying the Euclidean algorithm to the ratio $s/r$ where $s,r \in \mathbb{N}$ and $r<s$. Indeed, if
\begin{align}
s/ r  &  =q_{1}+u/r\nonumber\\
r/ u  &  =q_{2}+v/u\nonumber\\
&  \vdots\label{EA}\nonumber\\
x/ y  &  =q_{n}+0,\nonumber
\end{align}
then 
\begin{equation}s=[q_1,q_2,\ldots,q_n]\quad \text{and}\quad r=[q_2,q_3,\ldots,q_n].\nonumber
\end{equation}
Smith proves that if $p \equiv 1 \bmod4$, then there exists an integer $n$ and a unique $x_0 \in \mathbb{N}$ satisfying $1 < x_0 \leq (p-1)/2$ such that 
\begin{equation}
p=[q_1,q_2,\ldots,q_{n-1},q_n,q_n,q_{n-1},\ldots,q_{2},q_{1}] \nonumber
\end{equation}
and
\begin{equation}
x_0=[q_2,q_3\ldots,q_{n-1},q_n,q_n,q_{n-1},\ldots,q_{2},q_{1}]. \label{r}
\end{equation}
Notice that this continuant representation of $p$ is \emph{palindromic}. The integer $x_0$ is the smallest least residue of $x^2\equiv -1 \bmod{p}$. Furthermore, using elementary properties of determinants, Smith shows 
\begin{align}
p  &  =[q_{1},q_{2},\ldots q_{n-1},q_{n}][q_{n},q_{n-1},\ldots,q_{2},q_{1}]\nonumber\\
&  +[q_{1},q_{2},\ldots,q_{n-1}][q_{n-1},\ldots,q_{2},q_{1}%
]\label{Sum of Squares}\\
&  =[q_{1},q_{2},\ldots,q_{n}]^{2}+[q_{1},q_{2},\ldots,q_{n-1}]^{2}.\nonumber
\end{align}
Consequently, determining the least residue solution $x_0$ of $x^2 \equiv -1 \bmod{p}$ allows us, through the algorithm in (\ref{Sum of Squares}), to explicitly find positive integers $a,b$ satisfying $p=a^{2}+b^{2}$. To illustrate Smith's argument, consider $p=13$; in this case $x_0=5$ as can be seen from the following calculations using the Euclidean Algorithm, 
\[
\frac{13}{5}=\mathbf{2}+\frac{3}{5},\frac{5}%
{3}=\mathbf{1}+\frac{2}{3},\frac{3}{2}=\mathbf{1}+\frac{1}{2},\frac{2}%
{1}=\mathbf{2}+0,\]
so $p=[2,1,1,2]$ and $x_0=5=[1,1,2]$. It is clear that $x_0^2 \equiv -1 \bmod 13$. 
Moreover, employing (\ref{Sum of Squares}), we see that
\[
13   =[2,1,1,2]=[2,1][1,2]+[2][2] =[2,1]+[2]^{2}=3^{3}+2^{2}.
\]
As we will see in the next section, $x_0=5$ is called the \textit{St\o rmer number} for the prime $p=13$.

The key point in most proofs of Fermat's Two Squares Theorem is to first identify the positive integer $x_0$ satisfying the two conditions $1 < x_0 \leq (p-1)/2$ and $x_0^2 \equiv -1 \bmod p$. In practice, however, when the prime $p \equiv 1 \bmod 4$ is a large prime, it is difficult to find $x_0$; see Tables 1 and 2 below where it seems that there is considerable randomness in the list of St\o rmer numbers. In this paper, we obtain necessary and sufficient conditions for testing when $x_0 \in \mathbb{N}$ is a Størmer number and, in the case that it is, we show how to find the associated unique prime number $p$; see Theorem \ref{Main Theorem} below. 

The contents of this paper are as follows. Section \ref{Stormer bio} gives a brief biographical sketch of Carl St\o rmer; St\o rmer was a well-known Norwegian number theorist who was also recognized for his ground-breaking work in astronomy and his research on the aurora borealis. In Section \ref{Stormer Numbers}, we define St\o rmer numbers and prove several basic facts about them, including our characterization result (see Theorem \ref{Main Theorem}) of St\o rmer numbers. In Section \ref{Density}, we give a heuristic/probabilistic `proof' of the natural density of the St\o rmer numbers; it was conjectured by Everest and Harman \cite{Everest-Harman} in 2008 that this natural density is $\ln (2)$. Lastly, in Section \ref{Applications}, we discuss St\o rmer's work on approximating $\pi$ using the MacLaurin/Gregory series for $\arctan(1/x)$, that is,
\[
t_{x}:= \arctan \left(\frac{1}{x}\right)=\sum_{k=0}^\infty (-1)^k\frac{1}{(2k+1)x^{2k+1}}\quad(|x|\geq1). 
\]
Of course, it is well known that $t_1=\pi/4$. St\o rmer (see \cite{Stormer1896, Stormer1899}) shows that if $x \in \mathbb{N}$ is a \textit{non}-St\o rmer number, then $t_x$ can be written as a unique finite linear combination of $t_n$'s, where each $n$ is a St\o rmer number (see Theorem \ref{Stormer's Theorem}). We illustrate St\o rmer's method by considering several examples. 
\section{A Biographical Sketch of Carl St\o rmer}\label{Stormer bio}%

%TCIMACRO{\FRAME{fhFU}{2.0072in}{2.8236in}{0pt}{\Qcb{Carl St\o rmer}}%
%{}{stormer.jpg}{\special{ language "Scientific Word";  type "GRAPHIC";
%maintain-aspect-ratio TRUE;  display "USEDEF";  valid_file "F";
%width 2.0072in;  height 2.8236in;  depth 0pt;  original-width 1.9692in;
%original-height 2.7812in;  cropleft "0";  croptop "1";  cropright "1";
%cropbottom "0";  filename '../Stormer.jpg';file-properties "NPEU";}} }%
%BeginExpansion
\begin{figure}[h]%
\centering
\includegraphics[
height=2.8236in,
width=2.0072in
]{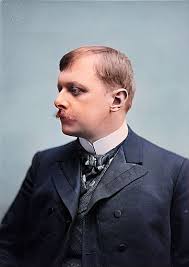}
\caption{Carl St\o rmer}%
\end{figure}
%EndExpansion
\ \ \ \ \ \ \  

Carl St\o rmer (1874-1957) was a preeminent Norwegian mathematician, botanist,
and astrophysicist. For his lifetime contributions to mathematics and astronomy,
particularly his study of the aurora borealis and the motion of charged particles in magnetic fields, he was bestowed several honors. A crater on the far side of the Moon is named after him. He was elected to
several Scandinavian academies, the Royal Society of London and the Paris
Academy of Sciences. He was given honorary degrees by the universities of
Oxford, Copenhagen and the Sorbonne. The Paris Academy of Sciences awarded him
their Janssen Medal in 1922. He was invited to give a one-hour lecture, on his
scientific  research on the aurora borealis, at the International Congress
of Mathematicians in Toronto in 1924 and, in 1936, he was president of the
International Congress of Mathematicians which was held in Oslo. 

St\o rmer entered the University of Christiania (Oslo's previous name) in 1892
and obtained his candidates' degree (similar to a Ph.D. degree) in 1897. By
this time, St\o rmer had written several mathematical articles, primarily in
number theory, as well as a number of short notes in botany. From 1898-1900,
St\o rmer studied at the Sorbonne under mathematical giants \'{E}mile Picard,
Henri Poincar\'{e}, Camille Jordan, Gaston Darboux, and \'{E}douard Goursat.
His output of mathematical papers continued with twelve papers on series,
number theory, and the theory of functions between 1896 and 1902. He worked
with Sylow and Holst to produce two volumes in 1902 to celebrate the centenary
of Niels Henrik Abel's (1802-1829) birth. 

Between 1896-1899, St\o rmer published two papers on connecting St\o rmer numbers to Gregory numbers, the subject of this paper. In 1897, he proved that if $P$ is a finite set of prime numbers, then there are only finitely many consecutive integers having only the numbers from $P$ as their prime factors. Moreover, St\o rmer develops an algorithm for finding all such pairs of consecutive integers. To prove his result, he reduced the problem to solving a finite number of Pell's equations. This paper \cite{Stormer1897} was highly praised by the number theorist Louis Mordell. 

In 1903, St\o rmer was appointed professor of pure mathematics at the University of Oslo, a position he held until 1946. In 1934, St\o rmer published a paper in the Norwegian Mathematical Society entitled `Ramanujan -- A Remarkable Genius'. At the insistence of fellow Norwegian Atle Selberg, a leading authority in number theory, this paper was reprinted in the same journal in 1989. Selberg was a teenager in 1934 and remarked that St\o rmer's paper on Ramanujan contributed to arousing his interest in mathematics.  

During his time at the
University of Oslo, St\o rmer's interests were also drawn to another
scientific area: the study of the aurora borealis. A colleague of St\o rmer,
Kristian Berkeland, had put forward a theory in 1896 that auroras were caused
by electrons emitted by the sun which interacted with the earth's magnetic
field. Poincar\'{e} had, in the same year, solved the differential equations
resulting from the motion of a charged particle in the field of a single pole.
This was not the situation for an aurora since the magnetic field of the earth
is a dipole. St\o rmer specifically attacked this problem. His papers on this
subject constitute a major part of his life's work. St\o rmer wrote two books,
\textit{From the Depths of Space to the Heart of the Atom}, which was
translated into five languages, and \textit{The Polar Aurora} which contains
both his experimental work on aurorae and his mathematical attempts to model
them. The 2013 book \textit{Carl St\o rmer: Aurora Pioneer}, by Egeland and
Burke \cite{Egeland and Burke} is a comprehensive account of St\o rmer's
contributions to the study of auroras.

We end this sketch with a rather humorous note about Carl St\o rmer. He was an avid
photographer which he put to use in his study of auroras. He purchased a
miniature \textit{spy} camera that he carefully concealed in his jacket. He
would walk the streets of Oslo and take photographs of people. He
published two works of his `snapshots of famous people' in 1942 and 1943. Long
afterwards, nearing the age of 70, these photographs formed the subject of a
major exhibition in Oslo. Many of these photographs exist on the internet and
they capture a fascinating look at life in Oslo during the early part of the 20th century. 

\section{St\o rmer Numbers}\label{Stormer Numbers}

\begin{definition}
\label{Definition of Stormer number}Suppose $p=4n+1$ is a prime number for some positive integer $n$. We
call a positive integer $x_{0}$ satisfying the two conditions
\begin{enumerate}
\item[(i)] $1<x_{0}\leq\dfrac{p-1}{2}$ (or $3<2x_{0}+1\leq p$)

\item[(ii)] $x_{0}^{2}\equiv-1\bmod p\medskip$ 
\end{enumerate} 
the St\o rmer number for $p$ and we write $S(p)=x_{0}$. If $\mathbb{P}_{4n+1}$ denotes the set of all prime numbers $p \equiv 1 \bmod4$, we call the function $S:\mathbb{P}%
_{4n+1}\rightarrow\mathbb{N}$, defined by $S(p)=x_{0}$, the St\o rmer function. We denote the set of St\o rmer numbers by $\mathbb{S}$ and its complement in $\mathbb{N}$ by $\mathbb{S}^{\text{c}}.$
By convention, we take $S(2)=1$.
\end{definition}

\noindent From the example in the last section, note that $S(13)=5$.

%\begin{remark}
%Conway and Guy's definition of a St\o rmer number (see \cite[p. 245]{Conway-Guy}) is slightly different from ours; they define a St\o rmer number to be a positive integer $n$ for which the largest prime factor $p$ of $n^2+1$ is at least $2n$ while our definition requires the largest prime factor to be at least $2n+1$. Their definition includes $n=1$ as a St\o rmer number; otherwise, the two definitions agree. 
%\end{remark}

St\o rmer numbers are listed as item A002314 in Sloane's website
\textit{Online Encyclopedia of Integer Sequences}. 

Table 1 below gives a list of a few ordered pairs $(p,S(p))$ in increasing order of $p$, where the prime
$p\equiv1\bmod 4$ and $S(p)$ is the corresponding St\o rmer number.%

\[%
\begin{tabular}
[c]{|l|l|l|l|l|l|l|}\hline
$(2,1)$ & $(5,2)$ & $(13,5)$ & $(17,4)$ & $(29,12)$ & $(37,6)$ & $(41,9)$ \\\hline
$(53,23)$ & $(61,11)$ & $(73,27)$ & $(89,34)$ & $(97,22)$ & $(101,10)$ & $(109,33)$\\\hline
$(113,15)$ & $(137,37)$ & $(149,44)$ & $(157,28)$ & $(173,80)$ & $(181,19)$ & $(193,81)$ \\\hline
$(197,14)$ &
$(229,107)$ & $(233,89)$ & $(241,64)$ & $(257,16)$ & $(269,82)$ & $(277,60)$ \\\hline
$(281,53)$ & 
$(293,138)$ & $(313,25)$ & $(317,114)$ & $(337,148)$ & $(349,136)$ &
$(353,42)$ \\\hline
\end{tabular}
\
\]
\begin{center}
Table 1
\end{center}
For example, $S(157)=28$ and $S(353)=42$. From this table, it appears that the relationship between $S(p)$ and $p$ is chaotic. As discussed in the introduction, we develop necessary and sufficient conditions for when a given $x_0 \in \mathbb{N}$ is the St\o rmer number for some prime $p \in \mathbb{P}_{4n+1}$. 

We list the first few St\o rmer numbers in increasing order:
\[
\begin{tabular}
[c]{|l|l|l|l|l|l|l|l|l|l|l|l|l|l|l|}\hline
$1$ & $2$ & $4$ & $5$ & $6$ & $9$ & $10$ & $11$ & $12$ & $14$ & $15$ & $16$ &
$19$ & $20$ & $22$\\\hline
$23$ & $24$ & $25$ & $26$ & $27$ & $28$ & $29$ & $33$ & $34$ & $35$ & $36$ &
$37$ & $39$ & $40$ & $42$\\\hline
$44$ & $45$ & $48$ & $49$ & $51$ & $52$ & $53$ & $54$ & $56$ & $58$ & $59$ &
$60$ & $61$ & $62$ & $63$\\\hline
$64$ & $65$ & $66$ & $67$ & $69$ & $71$ & $74$ & $77$ & $78$ & $79$ & $80$ &
$81$ & $82$ & $84$ & $85$\\\hline
$86$ & $87$ & $88$ & $90$ & $92$ & $94$ & $95$ & $96$ & $97$ & $101$ & $102$ &
$103$ & $104$ & $106$ & $107$\\\hline
\end{tabular}
\]
\begin{center}
Table 2
\end{center}

\begin{remark}
For a prime $p,$ it is well known by Euler's Criterion (see \cite[Theorems
9-1 and 9-5]{Andrews}) that the quadratic congruence
\begin{equation}
x^{2}\equiv-1\bmod p \label{S-1}%
\end{equation}
has exactly two least residue solutions when $p\equiv1\bmod 4$ (and
no solutions when $p\equiv3\bmod 4$). When $p\equiv
1\bmod 4,$ one of these solutions lies in the interval
$(1,(p-1)/2]$ and the other solution $y_{0}=p-x_{0}$ lies in the interval
$((p-1)/2,p)$.  Moreover, it is straightforward to see that $x_{0}=(p-1)/2$ satisfies (\ref{S-1}) only when $p=5$
(with $x_{0}=2).$ Hence, for primes $p>5$, which are
congruent to $1\bmod 4,$ its St\o rmer number $x_{0}$ satisfies
the strict inequality $1<x_{0}<(p-1)/2.$
\end{remark}
Our first result is that St\o rmer's function is injective.

\begin{theorem}
\label{S is 1-1}The St\o rmer function $S:\mathbb{P}_{4n+1}\rightarrow
\mathbb{N}$, defined in Definition \ref{Definition of Stormer number}, is 1-1;
that is, if $x_{0}=S(p)$ for some prime number $p,$ then $x_{0}$ cannot be the
St\o rmer number of any other prime $p^{\prime}\in\mathbb{P}_{4n+1}.$
\begin{proof}
Suppose $S(p_{1})=S(p_{2})$ but  $p_{1}\neq p_{2}$ for $j=1,2$. Since $p_{j}\vert (x_{0}^{2}+1)$, we
see that $p_{1}p_{2}\vert (x_{0}^{2}+1).$ But as $2x_0\leq p_j$, we see
\begin{equation}
p_1p_2\geq (2x_0)^2=4x_0^2
\end{equation}
and, hence, for some $k\in\mathbb{N}, x_0^2+1=k p_1 p_2\geq 4 x_0^2$. However, $x_0^2+1<4x_0^2$. Consequently, the St\o rmer map is $1$-$1$.
\end{proof}
\end{theorem}

It follows from Theorem \ref{S is 1-1} that there are infinitely many St\o rmer numbers; i.e. $\vert \mathbb{S} \vert=\aleph_0$, where $\vert A \vert$ denotes the cardinality of a set $A$.

For $x_0\in \mathbb{N}$, suppose the prime factorization of $x_0^{2}+1$ is%
\begin{equation}
x_0^{2}+1=2^{r}p_{1}^{r_{1}}p_{2}^{r_{2}}\cdots p_{m}^{r_{m}},
\label{Prime Factorization}%
\end{equation}
where $r$ is a non-negative integer, each $r_{j}\in\mathbb{N},$ and each
$p_{j}$ is an odd prime number with $p_{1}<p_{2}<\cdots<p_{m-1}<p_{m}.$

\begin{theorem}
Suppose $n=x_0^{2}+1$ has the prime factorization given in
(\ref{Prime Factorization}). Then, for each $1\leq j\leq m,$ $p_{j}%
\equiv1\bmod 4.$

\begin{proof}
Suppose, to the contrary, that $p_{j}\equiv3\bmod 4$ for some
$1\leq j\leq m.$ Since $p_{j}\vert (x_0^{2}+1),$ we see that $x_0^{2}\equiv
-1\bmod p_{j}.$ Furthermore, $\gcd(x_0,p_{j})=1$ (otherwise $p_{j}\vert x_0$ and
$p_{j}\vert (x_0^{2}+1)$ so $p\vert 1$ which is not possible), we see by Fermat's
Little Theorem that
\begin{equation}
x_0^{p_{j}-1}\equiv1\bmod p_{j}. \label{Fermat's Little Theorem}%
\end{equation}
Then, since $(p_{j}-1)/2$ is odd, we see that
\[
x_0^{p_{j}-1}=(x_0^{2})^{(p_{j}-1)/2}=(-1)^{(p_{j}-1)/2}\equiv-1\bmod %
p_{j},
\]
which contradicts (\ref{Fermat's Little Theorem}).
\end{proof}
\end{theorem}

We are now in position to characterize which positive integers $x_0$ can be the St\o rmer number of some prime $p\equiv1\bmod 4$.

\begin{theorem}
\label{Main Theorem}Suppose $n=x_0^{2}+1$ has the prime factorization given in
(\ref{Prime Factorization}) where we assume the primes satisfy $p_{1}%
<p_{2}<\cdots<p_{m}$. Then there exists a prime $p$ for which $x_0$ is the
St\o rmer number if and only if $2x_0+1\leq p_{m}$. Moreover, if this condition
is met, then $p=p_{m}$ and $S(p_{m})=x_0.$ Otherwise, $x_0$ is not the St\o rmer
number for any prime $p\equiv1\bmod 4.$

\begin{proof}
Suppose $2x_0+1\leq p_{m};$ that is,
\[
x_0\leq\frac{p_{m}-1}{2}.
\]
Since $p_{m}\vert (x_0^{2}+1),$ we see that $x_0$ is the St\o rmer number for
$p=p_{m};$ i.e. $S(p_{m})=x_0.\medskip$\newline Conversely, suppose $x_0=S(p)$ for
some prime $p\equiv1\bmod 4.$ It follows from the prime
factorization of $x_0^{2}+1$ that $p=p_{j}$ for some $j\in\{1,2,\ldots m\}$. By definition, $2x_0+1\leq p_{j}$ and
$x_0^{2}\equiv-1\bmod p_{j}.$ Suppose, for the sake of
contradiction, that $j<m$. Then since $p_{m}>p_{j}\geq2x_0+1$ and
$x_0^{2}\equiv-1\bmod p,$ we see that $x_0$ is also the St\o rmer
number for $p_{m}.$ But this contradicts the fact that $S$ is 1-1. This forces
$p=p_{m}.$ \medskip\newline If $2x_0+1>p_{m},$ then $x_0>(p_{m}-1)/2$ so, by
definition $x_0$ is not the St\o rmer number for $p_{m}$ nor, as the above
argument shows, can it be the St\o rmer number for any other prime
$p\equiv1\bmod 4.$
\end{proof}
\end{theorem}
\noindent \textbf{Examples}

\begin{enumerate}
\item Let $x_0=3.$ Since $x_0^{2}+1=3^{2}+1=2\cdot 5$ and $2x_0+1=7>5,$ $x_0=3$ is not
the St\o rmer number for any prime $p\equiv1\bmod 4.\smallskip$

\item Let $x_0=15.$ Since $x_0^{2}+1=226=2\cdot 113$ and $2x_0+1=31\leq113,$ it is
the case that $S(113)=15$ (see Table 1).$\smallskip$

\item Let $x_0=279.$ In this case, $x_0^{2}+1=2\cdot 38921$ and since
$2x_0+1=559\leq38921,$ we see that $x_0=279$ is the St\o rmer number for
$p=38921$.\smallskip

\item Suppose $p$ is a prime of the form $4n^{2}+1$ (for example, $5,17,$ or
$37).$ Then $x_0=2n$ is the St\o rmer number for $p=4n^{2}+1$ since
$2x_0+1=4n+1\leq4n^{2}+1.$
\end{enumerate}

\section{Comments on the Natural Density of St\o rmer Numbers}\label{Density}

The following table seems to suggest that, as $n\rightarrow\infty,$ the number $n$ of positive
integers which are St\o rmer numbers nears $70\%.$ 
\[%
\begin{tabular}
[c]{|l|l|}\hline
First $n$ positive integers & \# of St\o rmer numbers\\\hline
\multicolumn{1}{|c|}{$100$} & \multicolumn{1}{|c|}{$86$}\\\hline
\multicolumn{1}{|c|}{$1000$} & \multicolumn{1}{|c|}{$719$}\\\hline
\multicolumn{1}{|c|}{$10000$} & \multicolumn{1}{|c|}{$7101$}\\\hline
\multicolumn{1}{|c|}{$100000$} & \multicolumn{1}{|c|}{$70780$}\\\hline
\multicolumn{1}{|c|}{$1000000$} & \multicolumn{1}{|c|}{$704536$}\\\hline
\end{tabular}
\ \
\]

\begin{center}
Table 3
\end{center}
In fact, Everest and Harman \cite[Conjecture 1.5]{Everest-Harman} conjecture
\begin{equation}
\lim_{n \to \infty}\dfrac{\vert \{S(p) \mid  p \leq n\}\vert}{n} = \ln 2.\label{Conjecture}
\end{equation}

We now give a heuristic/probabilistic argument of (\ref{Conjecture}). We note that a rigorous proof still remains elusive. Let $x_0 \in \mathbb{N}$. If $x_0 = S(p)$ for some prime $p\equiv 1\bmod 4$, then $p \ge 2x_0+1$. Further,  since $x_0^2\equiv1\bmod p$ and $p\mid( x_0^2+1)$, we see that $p\le x_0^2+1$. Thus, for any prime $p \equiv 1\bmod4$, $S(p)$ can only be the St\o rmer number for $p$ if we have
\[
2x_0+1 \le p \le x_0^2+1.
\]
Assuming each of the integers from $1$ to $(p-1)/2$ is equally likely to be the St\o rmer number for $p$, we can informally think of the ``probability" that $x_0$ is a St\o rmer number for the prime $p$ as the reciprocal of $(p-1)/2$, that is, $2/(p-1)$. Then since $S(p)$ is injective, so that $x_0$ is a St\o rmer number for at most one prime, we can sum these all up to see that the ``probability" that $x_0$ is a St\o rmer number for any prime $p$ is
$$\sum_{\substack{2x_0+1\le p\le x_0^2+1 \\ p\equiv 1\bmod 4}}\frac2{p-1}.$$
This sum behaves like
$$\sum_{\substack{2x_0\le p\le x_0^2 \\ p\equiv 1\bmod 4}}\frac2p$$
since the effect of altering $p$ by 1 goes to zero in the limit. Since in the limit, there are the same number of primes of the form $p \equiv 1\bmod4$ and $p \equiv 3\bmod4$ (due to the Prime Number Theorem for arithmetic progressions), we write this as half the sum over all primes -- not just the primes of the form $p \equiv 1\bmod 4$ -- to obtain
$$\sum_{2x_0\le p\le x_0^2}\frac1p.$$
Then using Merten's estimate (see \cite[Lemma 4.10]{Bateman-Diamond}) we obtain
$$\sum_{p\le x_0}\frac1p \sim \ln\ln x_0,$$
we write this asymptotic as
$$\sum_{p\le x_0^2}\frac1p - \sum_{p< 2x_0}\frac1p \sim \ln(\ln x_0^2) - \ln(\ln 2x_0) \sim \ln\left(\frac{\ln x_0^2}{\ln 2x_0}\right) \sim \ln\left(\frac{2\ln x_0}{\ln x_0}\right) \sim \ln 2.$$
The gap in rigor in this argument is, of course, the introduction of probability and the assumption that $x_0$ has a $2/(p-1)$ chance of being the St\o rmer number for any given prime $p$.

\section{Pure Nostalgia: St\o rmer's Application of his Numbers}\label{Applications}
St\o rmer's interests in St\o rmer numbers originated from a problem that stretches back almost three millennia, specifically, the problem of accurately increasing the decimal digits of $\pi$. In the early 1900's, there were less than 600 digits of $\pi$ known. As of December 2025, the StorageReview Lab Team (see \cite{StorageReview2025} has computed more than 314 \textit{trillion} digits of $\pi$ known using A. Yee's \cite{Yee} y-cruncher software program. 

Throughout the history of calculating $\pi$ there have been, perhaps arguably, four major computational techniques used. These methods, chronologically, include Archimedes's approximation to $\pi$ using regular inscribed or circumscribed $n$-polygons, Machin and Machin-like formulas involving Gregory series, Gauss-Legendre methods, and most recently, Chudnovsky's algorithm. In this section, we focus on St\o rmer's work connecting Machin-like formualas and St\o rmer numbers. For further information on St\o rmer's work, we recommend \cite{Conway-Guy} and the original sources \cite{Stormer1896, Stormer1899}. Details on these other methods can be found in \cite{Bailey-Borwein} and \cite{BBB}. The mathematics behind St\o rmer's technique is fascinating and deserves attention not because of its application to approximating $\pi$ compared to today's superior methods but for its sheer mathematical elegance and beauty. 

In the $18^{\text{th}}$ century, due to the work of John Machin and Leonhard Euler, techniques were developed to approximate $\pi$ using \textit{Gregory numbers} $t_{x}$, defined by%
\[
t_{x}:=\arctan\left(  \frac{1}{x}\right)  =\sum_{k=0}^{\infty}(-1)^{k}\frac
{1}{(2k+1)x^{2k+1}} \qquad(\vert x\vert\geq 1).
\]
Of course, it is well-known that that $t_{1}=\pi/4$; this particular series is known as the M\={a}dhava-Gregory-Leibniz series. Convergence to $\pi/4$ is slow; indeed, roughly speaking, one must take 10 times more terms to add one extra decimal place of accuracy in the computation of $\pi$. Euler established the identity
\begin{equation}
t_{1}=5t_{7}+2t_{79/3}, \label{Euler}  
\end{equation} while
Machin proved that 
\begin{equation}
t_{1}=4t_{5}-t_{239} \label{Machin}
\end{equation}
which is equivalent to the well-known trigonometric 
identity%
\begin{equation}
\pi=16\arctan\left(  \frac{1}{5}\right)  -4\arctan\left(  \frac{1}%
{239}\right)  .
\end{equation}
Convergence is faster with Machin's formula; indeed, the first 100 terms of the series on the right-hand side of Machin's formula (\ref{Machin}) produces the first 140 digits of $\pi$. Compare this to Chudnovsky's algorithm which produces about $14$ new decimal digits of accuracy for every single term (or iteration) added to the Ramanujan-Sato series
\[
\dfrac{1}{\pi} = \dfrac{1}{12}\sum_{n=0}^\infty (-1)^n \dfrac{(6n)!(13591409+545140134n)}{(3n)!(n!)^3640320^{3n+\frac{3}{2}}}.
\]

In the late $19^{\text{th}}$ century, Carl St\o rmer noticed an important connection between Gaussian primes and identities involving Gregory numbers. Among his discoveries, in 1896, St\o rmer proved the identity%
\begin{equation}
t_{1}=44t_{57}+7t_{239}-12t_{682}+24t_{12943}.\label{Stormer identity}
\end{equation}
We note that, in 2002, the Japanese computer scientist Yasumasa Kanada and his team used this identity to help obtain the first $1.2411$ trillion digits of $\pi$ which, at the time, was the world record. 
 
Recall that a Gaussian integer $z=a+bi$, where $a,b \neq 0$, is a Gaussian prime if $\left\|z\right\|:=a^2+b^2$ is prime. For example, $i\pm 1$ and $3i\pm 2$ are Gaussian primes but $3\pm i$ is not. We note that with regard to the Gaussian primes, $\pm 1$, and $\pm i$, $\mathbb{Z}[i]$ forms a unique factorization domain. If $z=a+bi$, where $a>0$ and $\theta = \arg(z)$, where $\arg(z) \in (-\pi/2, \pi/2)$, then $\theta = \arctan(b/a)$. In particular, for $n \in \mathbb{N}$, we see that
\begin{equation}
\arg(n+i)=\arctan\left(\frac{1}{n}\right)=t_n.\label{arg 1}
\end{equation}
Moreover, since
\begin{equation}
\arg(z_1z_2)=\arg(z_1)+\arg(z_2),
\end{equation} we see, for $c>0$ and $k \in \mathbb{N}$, that
\begin{align}
\arg(cz)&=\arg(z)\label{arg 2}\\
\arg(z^k)&=k\arg(z).\label{powers of arg}
\end{align}
Now suppose that $t_n=t_{n_1}+t_{n_2}$, where $n,n_1, n_2 \in \mathbb{N}$. Then, from (\ref{arg 1}),
\begin{equation}
\arg((n_1+i)(n_2+i))=t_{n_1}+t_{n_2}=t_n. \label{arg 3}
\end{equation}
With an abuse of notation, but for better readability, we shall write (\ref{arg 3}) as
\begin{equation}
n+i \equiv (n_1+i)(n_2+i).\label{Stormer identity 2}
\end{equation}
This is the same notation adopted by Conway and Guy \cite{Conway-Guy}.

\begin{theorem}[Størmer]\label{Stormer's Theorem}
If $\gcd(a,b)=1$, then $t_{a/b}$ can be uniquely expressed as a finite linear combination of $t_n$'s where each $n$ is a Størmer number.
\end{theorem}

It is worth noting that many sources, including Conway and Guy \cite{Conway-Guy}, refer to St\o rmer's Theorem without an exact reference. In light of this, we now provide a proof of St\o rmer's Theorem which we believe to be new.

\begin{proof}
        Suppose $\gcd(a,b)=1$. Recall $t_{a/b}=\arctan\frac{b}{a}=\arg(a+bi)$. Then
        \[
        t_{a/b}=\pm c_1 t_{n_1}\pm c_2 t_{n_2}
        \pm\cdots\pm c_m t_{n_m}
        \]
        if and only if
        \[
        \arg(a+bi)=\arg((n_1\pm i)^{c_1}(n_2\pm i)^{c_2}\cdots(n_m\pm i)^{c_m}),
        \]
        or equivalently $(a+bi)\equiv(n_1\pm i)^{c_1}(n_2\pm i)^{c_2}\cdots(n_m\pm i)^{c_m}$, where congruence is taken as equivalence of arguments. It therefore suffices to show such a set of $n_j$'s and $c_j$'s exists and is unique.
        
          We first show the existence of such a set by induction over the norm of Gaussian integers. The Gaussian integers with coprime coefficients of least norm are $1\pm i$ and its associates, that is $1\pm i$ multiplied by a Gaussian unit. Recall that every Gaussian unit may be written as a power of $(1+i)$ or $(1-i)$. As $1$ is a Størmer number, this establishes the base case.
        
        As argument is unaffected by positive integer scalars, we may always assume the coefficients of a Gaussian integer are coprime. Now, suppose the theorem holds for all Gaussian integers of lesser norm than $(a+bi)$. Let $p$ be the largest prime factor of $a^2+b^2$ and let $x=S(p)$ (we note all prime divisors of $a^2+b^2$ are congruent to $1\pmod 4$ since $\gcd(a,b)=1$).
        Then $a^2\equiv -b^2\equiv b^2 x^2\pmod p$. It follows that $a=kp\pm bx$ for some $k\in\mathbb{Z}$ and $a\equiv \pm bx\pmod p$. Then 
        \[
            (a+bi)(x\mp i)=(kpx\pm b(x^2+1)\mp ikp).
        \]
        As $p\mid x^2+1$, $p\mid (kpx\pm b(x^2+1)\mp ikp)$ as well. Now consider
        \[
            \|a+bi\|-\left\|\frac{(kpx\pm b(x^2+1)\mp ikp)}{p}\right\|=\frac{(p^2-(x^2+1))((pk\pm bx)^2+ b^2)}{p^2}.
        \]
        As $p>x$, it follows that the numerator, and therefore the difference is positive for all primes $p$.
    
       Hence, by this process we yield a Gaussian integer of strictly lesser norm to which we may apply the induction hypothesis, thus proving the existence of Størmer decompositions. 

       We now consider the issue of uniqueness. Let $p$ be any prime factor of $a^2+b^2$ and let $x=S(p)$. Suppose there exists $y\in\mathbb{S}$, such that $p\mid y^2+1$, but $y\neq S(p)$. Let $q=S^{-1}(y)$. Then $q>p$ so $q\nmid (a^2+b^2)$. As a result, $q\nmid (a+bi)(x\pm i)$, but $q\mid\|(a+bi)(y\pm i)\|$. To eliminate $q$ we have two options:
        \begin{enumerate}
            \item Use $(y\mp i)$, which contributes $\arg(y+i)+\arg(y-i)=\arg(y^2+1)=0$ since $y^2+1$ is a positive real. This simply erases $(y\pm i)$ and contributes nothing to the decomposition.
            \item Choose $z\in\mathbb{S}$ with $q\mid z^2+1$ and $z\neq S(q)$, introducing a yet larger prime $S^{-1}(z)>q>p$.
        \end{enumerate}
        To obtain a different linear combination from the algorithm, we must therefore repeat the second option. But this produces a strictly increasing sequence of primes
        \[
        p < q < S^{-1}(z) < \ldots,
        \]
        which cannot terminate, contradicting the fact that each decomposition is finite. Hence $(x\pm i)=(S(p)\pm i)$ is the unique Størmer factor eliminating $p$. As the Størmer function is injective, uniqueness follows from unique factorization in the integers.
    \end{proof}
    
We note that in Machin's formula (\ref{Machin}), $239$ is not a St\o rmer number but both 1 and 5 are while, in St\o rmer's formula (\ref{Stormer identity}), none of the integers 57, 239, 682, 12943 are St\o rmer numbers. 

In the examples that follow, we use our notation to find identities involving Gregory numbers $t_n$. 

\begin{example}\label{Example 5.1}
Consider $3+i$. The prime factorization is $\|3+i\|^2=2\cdot 5$. Recalling $S(5)=2$ and noting $(2-i)\mid (3+i)$, we quickly see
\[
    (3+i)(2+i)=1+i.
\]
So $t_3= t_1-t_2$.
\end{example}

\begin{example}
    %The complex number $18-5i$ is of the form (II). To reduce it to a complex number of form (I), we solve the Diophantine equation $18d-5c=1$ to obtain $(c,d)=(7,2)$ which corresponds to $7+2i$. Since this latter number is not of the form (I), we solve $7d+2c=1$ to obtain $(c,d)=(4,-1)$ corresponding to $4-i$ which is of the form (I). 

    Consider the Gaussian integer $18+5i$ and note $\|18+5i\|^2=349$ a prime. We may calculate $S(349)=136$. In order to divide out the factor $349$, we must introduce the conjugate $18-5i$. As $(18-5i)\mid (136+i)$, we see
    \[
        (18+5i)(136+i)=349(7+2i).
    \]
    Repeating this process with $7+2i$, we see $\|7+2i\|^2=53$ a prime. We calculate $S(53)=23$ and note $(7-2i)\mid (23+i)$. So
    \[
        (18+5i)(136+i)(23+i)=349\cdot 53\cdot (3+i).
    \]
    Finally, $(3+i)=(1+i)(2-i)$ by Example \ref{Example 5.1}. Hence,
    \[
    t_{18/5}=t_1-t_2-t_{136}-t_{23}.
    \]
\end{example}

We note that others have studied St\o rmer numbers in their connection to the arctangent or cotangent functions. Since $\arctan(1/x)=\arccot(x)$, Lehmer \cite{Lehmer} studied identities of the form
\begin{equation*}
\arccot(a/b) = \arccot(n_0)-\arccot(n_1)+\arccot(n_2)-\cdots
\end{equation*}
where the integers $n_j$ are obtained by solving the recurrences
\[
a_j=n_jb_j+b_{j+1}\quad(0 \le b_{j+1}<b_j),\quad a_{j+1}=a_jn_j+b_j
\]
with initial conditions $a_0=a,b_0=b$.
Todd (see \cite{Todd} and the references cited therein) uses the notation $(n)=\arctan(n)$ instead of $t_n = \arctan(1/n)$ so, in our language, $(n)=t_{1/n}$. Todd defines, for $n \in \mathbb{N}$, ($n$) to be \textit{irreducible} if it cannot be written as a finite sum of the form
\[
\sum \alpha_j(n_j),
\]
where each $\alpha_j \in \mathbb{Z}$; otherwise, $(n)$ is \textit{reducible}. He proves the following theorem (see \cite[Theorem B]{Todd}) which is a slightly weaker version of Theorem \ref{Stormer's Theorem}. 
\begin{theorem}
For $n \in \mathbb{N}$, $n$ is a Størmer number if and only if $(n)$ is irreducible.
\end{theorem}

We remark that, using Todd's notation, Example \ref{Example 5.1} shows $(3)=3(1)-(2)$. 

In \cite[Theorem A]{Todd}, Todd also proves the following characterization theorem:
\begin{theorem}
For $n \in \mathbb{N}$, $n$ is a Størmer number if and only if all prime factors of $1+n^2$ occur among the prime factors of $1+m^2$ for $m=1,2,\cdots n-1$.
\end{theorem}

\begin{example}\label{Example 5.3}
The Slovene mathematician Jurij Vega held the world record for most calculated digits of $\pi$ from 1789 until 1841 using the following identity:
\[
t_1=2t_3+ t_7.
\]
From the factorization $\|7+i\|^2=2\cdot 5^2$ and the fact that $(2+i)\mid (7+i)$, we see
\[
    (7+i)(2-i)=5(3-i)\quad\text{or, equivalently,}\quad t_7+t_3=t_2.
\]
Recall, from Example \ref{Example 5.1}, that $t_3=t_1-t_2$. Adding these two identities we see $t_1=2t_3+t_7$.
\end{example}
%\begin{example}
   % The number $n=21$ is not a St\o rmer number by Theorem \ref{Main Theorem} since $21^2+1=2\cdot 13 \cdot 17$ and $2n+1=43 > 17$. The Gaussian Prime factorization of $21+i$ is,
%    \begin{equation}
%    21+i=(-i)(1+i)(3+2i)(4+i)\label{twenty one}.
%    \end{equation}
%    Both $1$ and $4$ are St\o rmer numbers so we only need to convert $-i$ and $3+2i$. A solution of the equation $3d+2c= -1$ is $(c,d)=(1,-1) = 1-i$. Multiplying both sides of (\ref{twenty one}) by $1-i$, we find that  
%    \[
  %  (1-i)(21+i)=(-i)(1+i)(5-i)(4+i).
 %   \]
 %   Next we use key identity (e) above, namely $(-i)=(1-i)^2$, to obtain
%    \[
%    (1-i)(21+i)=(1-i)^2(1+i)(5-i)(4+i).
 %   \]
 %   To isolate $21+i$, we multiply both sides by $1+i$ and obtain
  %  \[
 %   21+i=(1-i)^2(1+i)^2(5-i)(4+i)=(5-i)(4+i)
 %   \]
 %   and this yields the identity $t_{21}=t_4-t_5$.
%\end{example}

\begin{example}
    Since $70^2+1=13^2\cdot 29$, Theorem \ref{Main Theorem} says that $n=70$ is not a St\o rmer number. Note that $(2+5i),(2+3i)\mid (70+i)$. We calculate $S(13)=5$ and $S(29)=12$. Next, observe that $(2-3i)\mid (5-i)$ and $(2-5i)\mid (12-i)$. This yields
    \[
        (70+i)(5-i)^2(12-i)=9802(2-i).
    \]
    Hence, $t_{70}=-t_2+2t_5+t_{12}$.
    \end{example}

\begin{example}
We use St\o rmer's technique to establish Euler's identity in (\ref{Euler}). We begin by observing $\|79+3i\|^2=2\cdot5^5$; a calculation shows that $(1-i),(2-i)\mid (79+3i)$. Then
\[
(79+3i)(1+i)(2+i)^5=-6250.
\]
We note that $(1+i)^4\equiv -1$, so $t_{79/3}=3t_1-5t_2$. Hence, using Example \ref{Example 5.3},
\begin{align*}
    2t_{79/3}&=6t_1-10t_2\\
        &=t_1-5(t_1-2t_2)\\
        &=t_1-5t_7.
\end{align*}
\end{example}
\begin{example}
In our last example, we verify Machin's identity in (\ref{Machin}). The prime factorization of $\|239+i\|^2$ is given by
\[
239^2+1=2\cdot 13^4.
\]
Since $S(13)=5$ and $(5+i)\mid (239+i)$, we see
\[
(239+i)(5-i)^4=114244(1-i).
\]
and so $t_{239}=-t_1+4t_5$.
\end{example}
\vspace{.25in}
\noindent \textbf{Note:} Graeme Reinhart recently completed his undergraduate studies in mathematics at Baylor University.
\vspace{.1in}\newline
\noindent Author Contributions: Matthew Kroesche, Lance Littlejohn, and Graeme Reinhart equally contributed to the results in this manuscript and each of these authors helped in the final editing process. All authors agree to be accountable for all aspects of the work.
\vspace{.1in}\newline
Disclosure of interest: The authors have no competing interests to declare.
\vspace{.1in}\newline
Funding: No funding was received.

%\bibliographystyle{abbrv}
%\nocite{*}
%\bibliography{refs}

\end{document}